\documentclass[12pt]{amsart}
\usepackage{amsfonts, amsmath, amsthm}

\usepackage{graphicx}
\usepackage{epstopdf}
\usepackage{psfrag}

\newtheorem{pro}{Proposition}
\newtheorem{thm}[pro]{Theorem}

\newtheorem{cor}[pro]{Corollary}

\theoremstyle{definition}

\theoremstyle{remark}

\newcommand{\gen}{\mbox{\rm genus}}

\newcommand{\bdy}{\partial}

\newcommand{\thick}[1]{{\rm Thick}(#1)}
\newcommand{\thin}[1]{{\rm Thin}(#1)}

\title{Stabilizations of Heegaard splittings of sufficiently complicated 3-manifolds (Preliminary Report)} 
\date{\today}
\address{Pitzer College}
\email{bachman@pitzer.edu}
\author{David Bachman}
\begin{document}

\begin{abstract}
We construct several families of manifolds that have pairs of genus $g$ Heegaard splittings that must be stabilized roughly $g$ times to become equivalent. We also show that when two unstabilized, boundary-unstabilized Heegaard splittings are amalgamated by a ``sufficiently complicated" map, the resulting splitting is unstabilized. As a corollary, we produce a manifold that has distance one Heegaard splittings of arbitrarily high genus. Finally, we show
that in a 3-manifold formed by a sufficiently complicated gluing, a low genus, unstabilized Heegaard splitting can be expressed in a unique way as an amalgamation over the gluing surface.
\end{abstract}

\maketitle


\markright{SPLITTINGS OF SUFFICIENTLY COMPLICATED 3-MANIFOLDS}

\section{Introduction.}

The purpose of this report is to outline the main results of a series of forthcoming papers by the author. The first three of these \cite{dividing}, \cite{polyhedra}, \cite{NormalIndexN} lay the groundwork for ``topological index $n$" surfaces. The final paper in the series \cite{results} will replace the present one, and uses this technology to establish the results described below. 

Given a Heegaard surface $H$ in a 3-manifold, $M$, one can {\it stabilize} to obtain a splitting of higher genus by taking the connected sum of $H$ with the genus one splitting of $S^3$. Suppose $H_1$ and $H_2$ are Heegaard splittings of $M$, where $\rm{genus}(H_1) \ge \rm{genus}(H_2)$. It is a classical result of Reidemeister \cite{reidemeister} and Singer \cite{singer} from 1933 that as long as $H_1$ and $H_2$ induce the same partition of the components of $\partial M$, stabilizing $H_1$ some number of times produces a stabilization of $H_2$. Just one stabilization was proved to always be sufficient in large classes of 3-manifolds, including Seifert fibered spaces \cite{schultens:96}, genus two 3-manifolds \cite{rs:99}, and most graph manifolds \cite{Derby-Talbot2006} (see also \cite{sedgwick:97}). The lack of examples to the contrary has led to ``The Stabilization Conjecture": Any pair of Heegaard splittings requires at most one stabilization to become equivalent. (See Conjecture 7.4 in \cite{ScharlemannSurvey}.)

In this announcement we outline our program to construct families of counter-examples to the Stabilization Conjecture. This work, however, is broader than just the Stabilization Conjecture. The techniques developed to produce these examples also lead to the resolution of several other important questions in 3-manifold topology, as detailed below. 

This work was announced in December of 2007 at a Workshop on {\it Triangulations, Heegaard Splittings, and Hyperbolic Geometry}, at the American Institute of Mathematics. At the same conference another family of counter-examples to the Stabilization Conjecture was announced by Hass, Thompson, and Thurston, and their preprint has since appeared on the arxiv \cite{HTT}. Their proof uses mainly geometric techniques. Several months later Johnson posted a preprint on the arxiv \cite{johnson} that claims to be a combinatorial version of their paper. 

Our main family of counter-examples is described by the following theorem:

\begin{thm}
\label{t:MainFamily}
For each $n>5$ there is a closed 3-manifold which has a pair of inequivalent (unoriented) Heegaard splittings of genus $n$ which must be stabilized at least $n-3$ times to become equivalent. 
\end{thm}

The pairs of splittings obtained by Hass, Thompson, and Thurston are isotopic, as unoriented surfaces. When the manifold has boundary oriented and unoriented isotopy classes are the same. Johnson uses this fact to produce unoriented counter-examples. In fact, when there is boundary present our construction is easier as well, and is outlined in Section \ref{s:Program} below. This sketch is representative of the techniques used to establish all of the results claimed here.

\section{Higher Genus Gordon Conjecture}

Given an analogous result about connected sums (see \cite{GordonConjecture}), it is natural to conjecture that amalgamations of unstabilized Heegaard splittings are unstabilized. That is, if $M_1$ and $M_2$ are glued along a surface $F$ of non-zero genus, and $H_i$ is an unstabilized Heegaard surface in $M_i$, then $H_1$ and $H_2$ can be amalgamated in $M_1 \cup M_2$ to an unstabilized Heegaard surface. Unfortunately, even if the genus of $F$ is one, then Schultens and Weidmann have shown this to be false \cite{SchultensWeidmann}. It is perhaps surprising, then, that we are able to establish the following:

\begin{thm}
\label{t:HighGenusGordon}
Let $M_1$ and $M_2$ denote irreducible, orientable, anannular 3-manifolds with incompressible, homeomorphic boundary. Let $M$ be the 3-manifold obtained from $M_1$ and $M_2$ by gluing their boundaries by some ``sufficiently complicated" homeomorphism. Let $H_i$ be an unstabilized, boundary-unstabilized Heegaard surface in $M_i$ of low genus. Then the amalgamation of $H_1$ and $H_2$ in $M$ is unstabilized. 
\end{thm}

Just how low ``low genus" is depends on how complicated the gluing map between $M_1$ and $M_2$ is, unless $\bdy M_1 \cong T^2$. In this latter case the words ``low genus" may be removed from the statement of the theorem. 

As a corollary to this result we construct an example whose existence has been conjectured by Yoav Moriah: a non-minimal genus Heegaard splitting which has Hempel distance \cite{hempel:01} exactly one.  Moriah has called the search for such examples the ``nemesis of Heegaard splittings" \cite{moriah}. In fact, we go further and produce a single manifold that has such splittings of arbitrarily high genus:

\begin{cor}
There exists a closed 3-manifold that contains non-minimal genus, unstabilized Heegaard splittings which are not strongly irreducible, of arbitrarily high genus. 
\end{cor}

\begin{proof}
Let $M$ denote a 3-manifold with torus boundary, and strongly irreducible Heegaard splittings of arbitrarily high genus. Such an example has been constructed by Casson and Gordon. (See \cite{sedgwick:97}.)

Now let $M_1$ and $M_2$ be two copies of $M$, and let $H_g^i$ denote a genus $g$ strongly irreducible splitting in $M_i$. As $H_g^i$ is strongly irreducible, it is neither stabilized nor boundary-stabilized. Hence, if $M_1$ is glued to $M_2$ by a sufficiently complicated homeomorphism, it follows from Theorem \ref{t:HighGenusGordon} that the amalgamation of $H_g^1$ and $H_g^2$ is unstabilized. Finally, note that every amalgamation is weakly reducible. 
\end{proof}

Finally, we prove the following theorem, which is also analogous to a result about connected sums established by the author in \cite{GordonConjecture}. 

\begin{thm}
\label{t:HighGenusGordonIsotopy}
Let $M_1$ and $M_2$ denote irreducible, orientable, anannular 3-manifolds with incompressible, homeomorphic boundary.  Let $M$ be the 3-manifold obtained from $M_1$ and $M_2$ by gluing their boundaries by some ``sufficiently complicated" homeomorphism. Then every ``low genus" Heegaard splitting of $M$ has a \underline{unique} expression as the amalgamation of splittings of $M_1$ and $M_2$. 
\end{thm}

Again, the words ``low genus" in the statement of this theorem can be removed if $M_1$ and $M_2$ are being glued along a torus. 

Theorems \ref{t:HighGenusGordon} and \ref{t:HighGenusGordonIsotopy} together complete our picture of the ``low genus" Heegaard splittings of a 3-manifold constructed by gluing two manifolds $M_1$ and $M_2$ together by a ``sufficiently complicated" map.  Every unstabilized, low genus Heegard splitting of $M_1 \cup _\phi M_2$  is an amalgamation of ``component" unstabilized splittings of $M_1$ and $M_2$. Two such splittings are isotopic if and only if  their component splittings are isotopic. Hence, the unstabilized splittings of $M_1 \cup _\phi M_2$ are completely determined by the unstabilized splittings of $M_1$ and $M_2$.

\section{Normal and almost normal surfaces.}

A necessary step in the proofs of the above results is to put various classes of surfaces in some kind of ``normal form" with respect to a triangulation. We believe these results are of independent interest. 

Two common classes of surfaces that have proved useful in 3-manifold topology are the {\it incompressible} and {\it strongly irreducible} \cite{cg:87} ones. In \cite{crit} the author introduced a third class, called {\it critical surfaces}. These were used by the author in \cite{GordonConjecture} to prove a conjecture of Gordon. Each of these categories is defined by a combinatorial condition on the set of compressing disks for the surface. 

Following the work of Kneser \cite{kneser:29}, Haken proved that incompressible surfaces can be put into a nice form with respect to any triangulation \cite{haken:68}. Such surfaces were called {\it normal}, and intersect each tetrahedron in a collection of triangles and quadrilaterals. Rubinstein pushed this work further by proving that strongly irreducible Heegaard splittings can be isotoped to be {\it almost normal} (see also \cite{stocking:96}). Such surfaces look normal everywhere, with the exception of one piece in a single tetrahedron. This exceptional piece is either an octagon, or two normal disks connected by an unknotted tube. We call this piece the {\it almost normal exception}. 

In \cite{bachman:98} the author extended Rubinstein's result to strongly irreducible surfaces that have non-empty boundary (see also \cite{bachman:98erratum}). For this work, a third type of almost normal exception was introduced, which looks like two normal disks connected by a band along $\bdy M$. Similar results have been obtained by Wilson \cite{wilson} and Coward \cite{coward}. 

To prove the theorems listed above it was necessary to find a similar normal form for critical surfaces (possibly with non-empty boundary). The precise result is:

\begin{thm}
\label{t:CriticalNormal}
Given a fixed triangulation, any critical surface (with or without boundary) is isotopic to one which looks normal everywhere, except for exactly one 12-gon or two almost normal exceptions. 
\end{thm}

This theorem, as well as new proofs of the normalization results mentioned before it, will appear in \cite{NormalIndexN}.

\section{Topological Barrier Surfaces}

A technical result that will be of interest to the experts is that when two boundary components of a (possibly disconnected) 3-manifold are glued by a ``sufficiently complicated" map the gluing surface in the resulting 3-manifold acts as a barrier to a wide variety of other surfaces. A precise statement is:

\begin{thm}
\label{t:barrier}
 \cite{NormalIndexN} 
Let $M$ be a (possibly disconnected) 3-manifold with homeomorphic boundary components $F_1, F_2$. Let $M_\phi$ be the manifold obtained from $M$ by the gluing map $\phi:F_1 \to F_2$.  If $\phi$ is ``sufficiently complicated" then any incompressible, strongly irreducible, or critical surface in $M_\phi$ that has relatively low genus can be isotoped to lie entirely in $M$. 
\end{thm}

If $F_1 \cong T^2$ then the words ``relatively low genus" can be dropped from the statement of the above theorem. In all other cases the more complicated the map $\phi$ is, the higher the genus of the surface can be that we can guarantee will lie entirely in $M$. Theorem \ref{t:barrier} is the heart of the techniques discussed here. This is the main technical result used to obtain proofs of Theorems \ref{t:MainFamily} through \ref{t:HighGenusGordonIsotopy}. 

Theorem \ref{t:barrier} is a generalization of several previous results. Proofs of Theorem \ref{t:barrier} in the incompressible case have been given using both minimal surfaces and normal surfaces. Proofs in the strongly irreducible case have been given when $M$ is disconnected by Lackenby using (unstable) minimal surfaces \cite{lackenby:04} and Li using almost normal surfaces \cite{li}. There should certainly be a proof using minimal surfaces in the critical case, but it has only been conjectured that such surfaces can be isotoped to be minimal \cite{kyoto}. 

In lieu of this we use techniques from normal surface theory to establish the result. We prove in \cite{NormalIndexN} that if Theorem \ref{t:barrier} is false, then there are normal-like surfaces in $M$ whose boundaries on $F_1$ and $F_2$ are a bounded distance apart in the curve complex after gluing via $\phi$.  Theorem \ref{t:barrier} then follows from \cite{js:98} in the case where $F_1 \cong T^2$ and $M$ is disconnected, from \cite{self} when $F_1 \cong T^2$ and $M$ is connected,  and from \cite{jr:Half} when $F_1$ has higher genus.

\section{The Program}
\label{s:Program}

In this section we illustrate the techniques used to obtain Theorems \ref{t:MainFamily} through \ref{t:HighGenusGordonIsotopy}. Here we sketch our proof that there is a family of counter-examples to the Stabilization Conjecture that are splittings of a manifold with torus boundary. 

To begin, let $X$ be a compact, orientable, irreducible, annannular 3-manifold such that:
\begin{enumerate}
	\item One component of $\partial X$ is a torus $T$ which is incompressible in $X$.
	\item The other component of $\partial X$ is a high genus surface $F$, which is incompressible in $X$.
	\item There is a minimal genus Heegaard splitting $H_X$ of $X$ which separates $F$ from $T$, and has genus equal to $\gen(F)+1$.
\end{enumerate}

Such a manifold is not difficult to construct. We can obtain a new Heegaard splitting $G_X$ of $X$ by tubing $H_X$ to a copy of $T$ (pushed slightly into $X$). The resulting splitting does not separate the boundary components of $X$. It follows that $H_X$ and $G_X$ are not equivalent after {\it any} number of stabilizations in $X$. This observation is the key to the proof that the Heegaard splittings constructed presently give counterexamples to the Stabilization Conjecture.

Now let $Y$ be a 3-manifold with incompressible boundary, such that $\partial Y \cong F$. Let $H_Y$ denote a minimal genus Heegaard splitting of $Y$. Such a manifold can be constructed where $\gen(H_Y)=\gen(H_X)$. For each homeomorphism $\phi: \partial Y \to F$ let $M_{\phi}=X \cup _{\phi} Y$ denote the manifold obtained by gluing $X$ to $Y$ via $\phi$. 

It is well known that $H_X$ and $H_Y$ give rise to a Heegaard splitting $H_\phi$ of $M_{\phi}$, called their {\it amalgamation} \cite{schultens:93}. Similarly, there is an amalgamation of $G_X$ and $H_Y$, which we call $G_\phi$. 

\begin{thm}
\label{t:TorusBoundaryCounterExample}
If $\phi$ is ``sufficiently complicated" then one must stabilize $G_\phi$ at least $\gen(G_\phi)-4$ times to obtain a stabilization of $H_\phi$.
\end{thm}

This theorem is proved by using critical surfaces. To give a sketch of the proof, we must first describe how such surfaces arise. 

Loosely speaking, a {\it Generalized Heegaard Splitting} GHS $H$ of a 3-manifold $M$ is a pair of sets of orientable, connected surfaces, $\thick{H}$ and $\thin{H}$, such that each component $M'$ of the complement of $\thin{H}$ in $M$ contains exactly one element $H_+$ of $\thick{H}$, and $H_+$ is a Heegaard surface for $M'$. From the point of view of  \cite{st:94}, GHSs arise naturally from handle structures. Conversely, given a GHS one can always find an associated handle structure.

There are two natural ways to take a GHS and produce a ``simpler" one. These can both be described by the following procedure:
	\begin{enumerate}
		\item Beginning with a GHS, form an associated handle structure.
		\item Either (a) swap the order of attachment of a 2-handle that had followed a 1-handle, or (b) cancel a 1-handle and a 2-handle.
		\item Form the associated GHS of the new handle structure. 
	\end{enumerate}
In Case (a) we say the new GHS is obtained from $H$ by a {\it weak reduction}. In Case (b) we say it was obtained by {\it destabilization}. There is a complexity such that in either case the new GHS is ``smaller" than the original. If there are no weak reductions or destabilizations for a GHS, then each of its thick levels are strongly irreducible. Scharlemann and Thompson showed that this implies that its thin levels must then be incompressible. 

A {\it Sequence Of GHSs} (SOG) is defined to be any sequence where each element is obtained from its predecessor or successor by a weak reduction or destabilization. If an element $H$ of a SOG $\bf H$ is greater than both its predecessor and its successor then we say it is {\it maximal} in $\bf H$. 

Just as there are two ways to take a GHS and produce a smaller one, there are ways one can ``reduce" a SOG. If a SOG can not be reduced, then each thick level of each maximal GHS is either strongly irreducible or critical and each thin level of each maximal GHS is incompressible.

\begin{figure}
\psfrag{H}{$H_X$}
\psfrag{h}{$G_X$}
\psfrag{G}{$H_Y$}
\psfrag{F}{$F$}
\psfrag{T}{$T$}
\psfrag{1}{$H _\phi$}
\psfrag{2}{$G_\phi$}
\psfrag{*}{$H_*$}
\[\includegraphics[width=5 in]{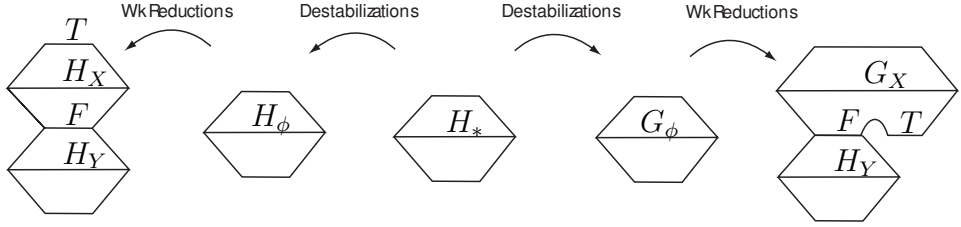}\]
\caption{The initial SOG $\bf H$. Note that $H_X$ separates $F$ and $T$ in $X$, while $G_X$ does not.}
\label{f:NSFInitialSOG}
\end{figure}

We are now prepared to sketch the proof of Theorem \ref{t:TorusBoundaryCounterExample}. First, observe that the basic set-up of the hypotheses gives us the SOG $\Gamma$ of $X \cup _\phi Y$ pictured in Figure \ref{f:NSFInitialSOG}. If there is no way to reduce the SOG $\Gamma$ then it would follow that the minimal genus common stabilization $H_*$ of $H_\phi$ and $G_\phi$ is critical. Since no Heegaard surface in $X \cup _\phi Y$ can be isotoped to be disjoint from $F$, it now follows from Theorem \ref{t:barrier} that the genus of $H_*$ must be very high. Hence, if Theorem \ref{t:TorusBoundaryCounterExample} is false then there must be some way to reduce the SOG $\Gamma$.

We now assume $\Lambda$ is obtained from $\Gamma$ by a maximal sequence of complexity reducing moves. It follows from Theorem \ref{t:barrier} that the surface $F$ must be an element of $\thin{\Lambda^i}$ for each GHS $\Lambda^i \in \Lambda$.

\begin{figure}
\psfrag{F}{$F$}
\psfrag{T}{$T$}
\[\includegraphics[width=5 in]{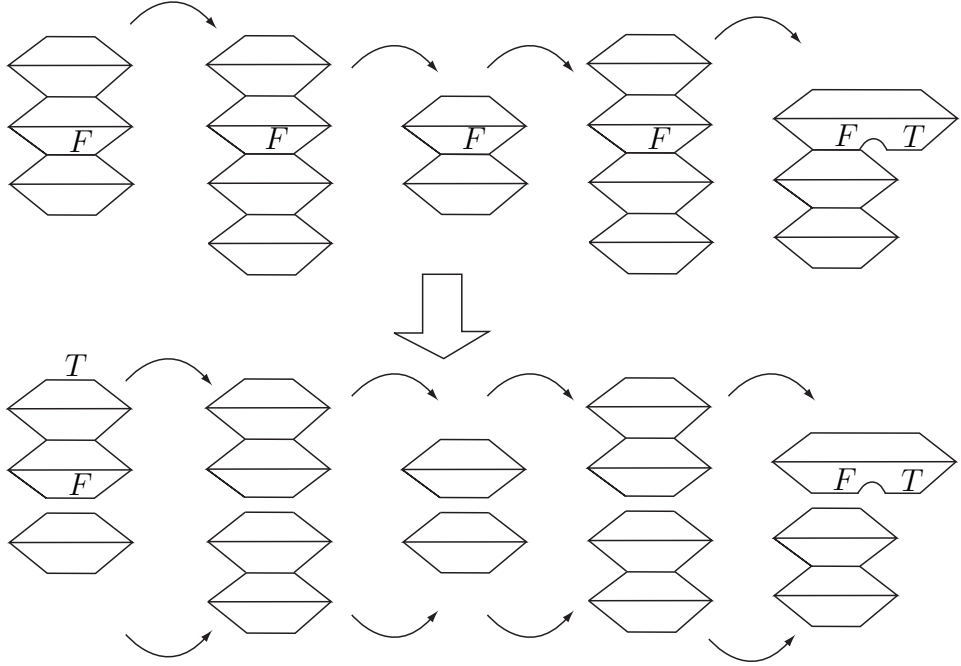}\]
\caption{An SOG with $F$ as a unique thin level of each GHS can be ``broken" at $F$ to obtain SOGs of $X$ and $Y$.}
\label{f:BreakF}
\end{figure}

If $F$ is a unique element of $\thin{\Lambda^i}$, for each element $\Lambda^i \in \Lambda$, then the SOG $\Lambda$ can be ``broken" at $F$ to obtain separate SOGs of $X$ and $Y$. See Figure \ref{f:BreakF}. The first GHS of the SOG of $X$ thus obtained comes from weak reductions and destabilizations of the Heegaard surface $H_X$. As such it separates the surfaces $F$ and $T$. The final GHS of this SOG is descended from $G_X$, and therefore does not separate $F$ and $T$. We now have a contradiction, as there can be no SOG which interpolates between GHSs that partition the boundary of $X$ differently. 

We conclude that there must be an element $\Lambda^*$ of $\Lambda$ for which the surface $F$ appears as two elements of $\thin{\Lambda^*}$. Hence, there is an element of $\thick{\Lambda^*}$ which lies between two copies of $F$, and is thus a Heegaard surface in a submanifold homeomorphic to $F \times I$. A classification Theorem of Scharlemann and Thompson \cite{st:93} tells us that this surface must be two copies of $F$ connected by a tube. (The classification includes one more possibility which is ruled out.)

Each GHS has a well-defined {\it genus}, which is the sum of the genera of its thick surfaces minus the sum of the genera of its thin surfaces. The above argument leads us to the inequality
\[genus(\Lambda^*) \ge genus(X)+genus(Y)\]
But, an analysis of the complexity reducing moves used to obtain $\Lambda$ now shows that there must be a GHS of $\Gamma$ whose genus is at least that of $\Lambda^*$. Finally, the genus of the surface $H_*$ is the largest genus of all the GHSs in $\Gamma$. Putting these inequalities together gives us $\gen(H_*) \ge \gen(X)+\gen(Y)$. But the genus of $G_\phi$ is $1+\gen(X)+\gen(Y)-\gen(F)$. The difference between these is the number of times we must stabilize $G_\phi$ in order to obtain a stabilization of $H_\phi$. This number is $\gen(F)-1$, which is also equal to $\gen(G_\phi)-4$. 
	

\end{document}